
\input mssymb


\def\0{\text{\bf 0}}
\def\Rn{{\Bbb R}^n}
\def\R2n{{\Bbb R}^{2n}}
\def\RnZ{(\Rn)^\Z}

\def\T1{{\Bbb S}^1}
\def\Tn{{\Bbb T}^n}
\def\TTn{T^*\Tn}

\def\Wo_{\underline{WO}_{p,q}}
\def\XmdN{X^*}

\def\xx{X}
\def\Z{{\Bbb Z}}
\def\Zn{{\Bbb Z}^n}
\def\zz{{\vec z}}

\def\F{\tilde F}

\def\pr{pr}     
\def\Tm{{\tau}_\m}


\def\WmdN{W}
\def\Hpp{H_{\p\p}}

\def\m{\vec{m}}
\def\q{\vec{ q}}
\def\qq{\overline \q}
\def\Q{\vec{ Q}}

\def\p{\vec{ p}}
\def\P{\vec{ P}}

\def\vv{\vec v}

\def\z{\vec z}


\def\abs#1{\left\vert #1\right\vert}  
\def\norm#1{\left\Vert #1\right\Vert}   


\def\del{\partial}  
\def\dd#1#2{{\partial #1\over \partial #2}}  

\def\t#1{\tilde #1}

\def\gf#1#2{ S_{#1}(\q_{#1},\q_{#2})} 

\def\qbox#1{\quad\hbox{#1}\quad}

\def\TheDate{\number\day~\ifcase\month\or
January\or February\or March\or April\or May\or June\or July\or
August\or September\or October\or November\or December\fi, \number\year}

\newcount\formulano
\newcount\sectionno
\newcount\theoremno
\newcount\figureno
\sectionno=0

\def\advsection{\global\theoremno=0\global\formulano=0\global\advance\sectionno
by1}
\def\secno{\bf \the\sectionno .}
\def\Secno#1{\global\advance\sectionno by #1 {\rm \the\sectionno }\global
\advance\sectionno by -#1}

\def\formno#1{\global\advance\formulano
by #1 {\rm(\the\sectionno .\the\formulano)\ }\ifnum #1 
<0\global\advance\formulano by -#1\fi}
\def\numberform{\leqno\formno 1}
\def\figno#1{\global\advance\figureno by #1{\bf \the\figureno}}
\def\thmno#1{\global\advance\theoremno by #1
{\bf \the\sectionno .\the\theoremno \  }\ifnum #1<0\global\advance\theoremno by 
-#1\fi}

\def\b{{\vec b}}

\def\RR{{\Bbb R}}

\def\NN{{\Bbb N}}
\def\ZZ{{\Bbb Z}}

\def\TT{{\Bbb T}}

\magnification=\magstep1
\font \authfont               = cmr10 scaled\magstep4
\font \fivesans               = cmss10 at 5pt
\font \headfont               = cmbx12 scaled\magstep4
\font \markfont               = cmr10 scaled\magstep1
\font \ninebf                 = cmbx9
\font \ninei                  = cmmi9
\font \nineit                 = cmti9
\font \ninerm                 = cmr9
\font \ninesans               = cmss10 at 9pt
\font \ninesl                 = cmsl9
\font \ninesy                 = cmsy9
\font \ninett                 = cmtt9
\font \sevensans              = cmss10 at 7pt
\font \sixbf                  = cmbx6
\font \sixi                   = cmmi6
\font \sixrm                  = cmr6
\font \sixsans                = cmss10 at 6pt
\font \sixsy                  = cmsy6
\font \smallescriptfont       = cmr5 at 7pt
\font \smallescriptscriptfont = cmr5
\font \smalletextfont         = cmr5 at 10pt
\font \subhfont               = cmr10 scaled\magstep4
\font \tafonts                = cmbx7  scaled\magstep2
\font \tafontss               = cmbx5  scaled\magstep2
\font \tafontt                = cmbx10 scaled\magstep2
\font \tams                   = cmmib10
\font \tamss                  = cmmib10 scaled 700
\font \tamt                   = cmmib10 scaled\magstep2
\font \tass                   = cmsy7  scaled\magstep2
\font \tasss                  = cmsy5  scaled\magstep2
\font \tast                   = cmsy10 scaled\magstep2
\font \tasys                  = cmex10 scaled\magstep1
\font \tasyt                  = cmex10 scaled\magstep2
\font \tbfonts                = cmbx7  scaled\magstep1
\font \tbfontss               = cmbx5  scaled\magstep1
\font \tbfontt                = cmbx10 scaled\magstep1
\font \tbms                   = cmmib10 scaled 833
\font \tbmss                  = cmmib10 scaled 600
\font \tbmt                   = cmmib10 scaled\magstep1
\font \tbss                   = cmsy7  scaled\magstep1
\font \tbsss                  = cmsy5  scaled\magstep1
\font \tbst                   = cmsy10 scaled\magstep1
\font \tenbfne                = cmb10
\font \tensans                = cmss10
\font \tpfonts                = cmbx7  scaled\magstep3
\font \tpfontss               = cmbx5  scaled\magstep3
\font \tpfontt                = cmbx10 scaled\magstep3
\font \tpmt                   = cmmib10 scaled\magstep3
\font \tpss                   = cmsy7  scaled\magstep3
\font \tpsss                  = cmsy5  scaled\magstep3
\font \tpst                   = cmsy10 scaled\magstep3
\font \tpsyt                  = cmex10 scaled\magstep3
\vsize=22.5true cm
\hsize=13.8true cm
\hfuzz=2pt
\tolerance=500
\abovedisplayskip=3 mm plus6pt minus 4pt
\belowdisplayskip=3 mm plus6pt minus 4pt
\abovedisplayshortskip=0mm plus6pt minus 2pt
\belowdisplayshortskip=2 mm plus4pt minus 4pt
\predisplaypenalty=0
\clubpenalty=10000
\widowpenalty=10000
\frenchspacing
\newdimen\oldparindent\oldparindent=1.5em
\parindent=1.5em
\skewchar\ninei='177 \skewchar\sixi='177
\skewchar\ninesy='60 \skewchar\sixsy='60
\hyphenchar\ninett=-1
\def\newline{\hfil\break}%
\catcode`@=11
\def\folio{\ifnum\pageno<\z@
\uppercase\expandafter{\romannumeral-\pageno}%
\else\number\pageno \fi}
\catcode`@=12 
  \mathchardef\Gamma="0100
  \mathchardef\Delta="0101
  \mathchardef\Theta="0102
  \mathchardef\Lambda="0103
  \mathchardef\Xi="0104
  \mathchardef\Pi="0105
  \mathchardef\Sigma="0106
  \mathchardef\Upsilon="0107
  \mathchardef\Phi="0108
  \mathchardef\Psi="0109
  \mathchardef\Omega="010A
  \mathchardef\bfGamma="0\the\bffam 00
  \mathchardef\bfDelta="0\the\bffam 01
  \mathchardef\bfTheta="0\the\bffam 02
  \mathchardef\bfLambda="0\the\bffam 03
  \mathchardef\bfXi="0\the\bffam 04
  \mathchardef\bfPi="0\the\bffam 05
  \mathchardef\bfSigma="0\the\bffam 06
  \mathchardef\bfUpsilon="0\the\bffam 07
  \mathchardef\bfPhi="0\the\bffam 08
  \mathchardef\bfPsi="0\the\bffam 09
  \mathchardef\bfOmega="0\the\bffam 0A

\def\sq{\hbox{\rlap{$\sqcap$}$\sqcup$}}

\def\utw{\smash{\rlap{\lower5pt\hbox{$\sim$}}}}
\def\udtw{\smash{\rlap{\lower6pt\hbox{$\approx$}}}}

\def\diameter{{\ifmmode\mathchoice
{\ooalign{\hfil\hbox{$\displaystyle/$}\hfil\crcr
{\hbox{$\displaystyle\mathchar"20D$}}}}
{\ooalign{\hfil\hbox{$\textstyle/$}\hfil\crcr
{\hbox{$\textstyle\mathchar"20D$}}}}
{\ooalign{\hfil\hbox{$\scriptstyle/$}\hfil\crcr
{\hbox{$\scriptstyle\mathchar"20D$}}}}
{\ooalign{\hfil\hbox{$\scriptscriptstyle/$}\hfil\crcr
{\hbox{$\scriptscriptstyle\mathchar"20D$}}}}
\else{\ooalign{\hfil/\hfil\crcr\mathhexbox20D}}%
\fi}}


\def\bbbc{{\mathchoice {\setbox0=\hbox{$\displaystyle\rm C$}\hbox{\hbox
to0pt{\kern0.4\wd0\vrule height0.9\ht0\hss}\box0}}
{\setbox0=\hbox{$\textstyle\rm C$}\hbox{\hbox
to0pt{\kern0.4\wd0\vrule height0.9\ht0\hss}\box0}}
{\setbox0=\hbox{$\scriptstyle\rm C$}\hbox{\hbox
to0pt{\kern0.4\wd0\vrule height0.9\ht0\hss}\box0}}
{\setbox0=\hbox{$\scriptscriptstyle\rm C$}\hbox{\hbox
to0pt{\kern0.4\wd0\vrule height0.9\ht0\hss}\box0}}}}
\def\bbbe{{\mathchoice {\setbox0=\hbox{\smalletextfont e}\hbox{\raise
0.1\ht0\hbox to0pt{\kern0.4\wd0\vrule width0.3pt height0.7\ht0\hss}\box0}}
{\setbox0=\hbox{\smalletextfont e}\hbox{\raise
0.1\ht0\hbox to0pt{\kern0.4\wd0\vrule width0.3pt height0.7\ht0\hss}\box0}}
{\setbox0=\hbox{\smallescriptfont e}\hbox{\raise
0.1\ht0\hbox to0pt{\kern0.5\wd0\vrule width0.2pt height0.7\ht0\hss}\box0}}
{\setbox0=\hbox{\smallescriptscriptfont e}\hbox{\raise
0.1\ht0\hbox to0pt{\kern0.4\wd0\vrule width0.2pt height0.7\ht0\hss}\box0}}}}
\def\bbbq{{\mathchoice {\setbox0=\hbox{$\displaystyle\rm Q$}\hbox{\raise
0.15\ht0\hbox to0pt{\kern0.4\wd0\vrule height0.8\ht0\hss}\box0}}
{\setbox0=\hbox{$\textstyle\rm Q$}\hbox{\raise
0.15\ht0\hbox to0pt{\kern0.4\wd0\vrule height0.8\ht0\hss}\box0}}
{\setbox0=\hbox{$\scriptstyle\rm Q$}\hbox{\raise
0.15\ht0\hbox to0pt{\kern0.4\wd0\vrule height0.7\ht0\hss}\box0}}
{\setbox0=\hbox{$\scriptscriptstyle\rm Q$}\hbox{\raise
0.15\ht0\hbox to0pt{\kern0.4\wd0\vrule height0.7\ht0\hss}\box0}}}}
\def\bbbt{{\mathchoice {\setbox0=\hbox{$\displaystyle\rm
T$}\hbox{\hbox to0pt{\kern0.3\wd0\vrule height0.9\ht0\hss}\box0}}
{\setbox0=\hbox{$\textstyle\rm T$}\hbox{\hbox
to0pt{\kern0.3\wd0\vrule height0.9\ht0\hss}\box0}}
{\setbox0=\hbox{$\scriptstyle\rm T$}\hbox{\hbox
to0pt{\kern0.3\wd0\vrule height0.9\ht0\hss}\box0}}
{\setbox0=\hbox{$\scriptscriptstyle\rm T$}\hbox{\hbox
to0pt{\kern0.3\wd0\vrule height0.9\ht0\hss}\box0}}}}
\def\bbbs{{\mathchoice
{\setbox0=\hbox{$\displaystyle     \rm S$}\hbox{\raise0.5\ht0\hbox
to0pt{\kern0.35\wd0\vrule height0.45\ht0\hss}\hbox
to0pt{\kern0.55\wd0\vrule height0.5\ht0\hss}\box0}}
{\setbox0=\hbox{$\textstyle        \rm S$}\hbox{\raise0.5\ht0\hbox
to0pt{\kern0.35\wd0\vrule height0.45\ht0\hss}\hbox
to0pt{\kern0.55\wd0\vrule height0.5\ht0\hss}\box0}}
{\setbox0=\hbox{$\scriptstyle      \rm S$}\hbox{\raise0.5\ht0\hbox
to0pt{\kern0.35\wd0\vrule height0.45\ht0\hss}\raise0.05\ht0\hbox
to0pt{\kern0.5\wd0\vrule height0.45\ht0\hss}\box0}}
{\setbox0=\hbox{$\scriptscriptstyle\rm S$}\hbox{\raise0.5\ht0\hbox
to0pt{\kern0.4\wd0\vrule height0.45\ht0\hss}\raise0.05\ht0\hbox
to0pt{\kern0.55\wd0\vrule height0.45\ht0\hss}\box0}}}}
\def\bbbz{{\mathchoice {\hbox{$\sans\textstyle Z\kern-0.4em Z$}}
{\hbox{$\sans\textstyle Z\kern-0.4em Z$}}
{\hbox{$\sans\scriptstyle Z\kern-0.3em Z$}}
{\hbox{$\sans\scriptscriptstyle Z\kern-0.2em Z$}}}}
\def\qed{\ifmmode\sq\else{\unskip\nobreak\hfil
\penalty50\hskip1em\null\nobreak\hfil\sq
\parfillskip=0pt\finalhyphendemerits=0\endgraf}\fi}
\newfam\sansfam
\textfont\sansfam=\tensans\scriptfont\sansfam=\sevensans
\scriptscriptfont\sansfam=\fivesans
\def\sans{\fam\sansfam\tensans}
\def\stackfigbox{\if
Y\FIG\global\setbox\figbox=\vbox{\unvbox\figbox\box1}%
\else\global\setbox\figbox=\vbox{\box1}\global\let\FIG=Y\fi}
\def\placefigure{\dimen0=\ht1\advance\dimen0by\dp1
\advance\dimen0by5\baselineskip
\advance\dimen0by0.4true cm
\ifdim\dimen0>\vsize\pageinsert\box1\vfill\endinsert
\else
\if Y\FIG\stackfigbox\else
\dimen0=\pagetotal\ifdim\dimen0<\pagegoal
\advance\dimen0by\ht1\advance\dimen0by\dp1\advance\dimen0by1.4true cm
\ifdim\dimen0>\pagegoal\stackfigbox
\else\box1\vskip4true mm\fi
\else\box1\vskip4true mm\fi\fi\fi}
%
\def\begfig#1cm#2\endfig{\par
\setbox1=\vbox{\dimen0=#1true cm\advance\dimen0
by1true cm\kern\dimen0#2}\placefigure}
\def\begdoublefig#1cm #2 #3 \enddoublefig{\begfig#1cm%
\vskip-.8333\baselineskip\line{\vtop{\hsize=0.46\hsize#2}\hfill
\vtop{\hsize=0.46\hsize#3}}\endfig}
\def\begfigsidebottom#1cm#2cm#3\endfigsidebottom{\dimen0=#2true cm
\ifdim\dimen0<0.4\hsize\message{Room for legend to narrow;
begfigsidebottom changed to begfig}\begfig#1cm#3\endfig\else
\par\def\figure##1##2{\vbox{\noindent\petit{\bf
Fig.\ts##1\unskip.\ }\ignorespaces ##2\par}}%
\dimen0=\hsize\advance\dimen0 by-.8true cm\advance\dimen0 by-#2true cm
\setbox1=\vbox{\hbox{\hbox to\dimen0{\vrule height#1true cm\hrulefill}%
\kern.8true cm\vbox{\hsize=#2true cm#3}}}\placefigure\fi}
\def\begfigsidetop#1cm#2cm#3\endfigsidetop{\dimen0=#2true cm
\ifdim\dimen0<0.4\hsize\message{Room for legend to narrow; begfigsidetop
changed to begfig}\begfig#1cm#3\endfig\else
\par\def\figure##1##2{\vbox{\noindent\petit{\bf
Fig.\ts##1\unskip.\ }\ignorespaces ##2\par}}%
\dimen0=\hsize\advance\dimen0 by-.8true cm\advance\dimen0 by-#2true cm
\setbox1=\vbox{\hbox{\hbox to\dimen0{\vrule height#1true cm\hrulefill}%
\kern.8true cm\vbox to#1true cm{\hsize=#2 true cm#3\vfill
}}}\placefigure\fi}
\def\figure#1#2{\vskip1true cm\setbox0=\vbox{\noindent\petit{\bf
Fig.\ts#1\unskip.\ }\ignorespaces #2\smallskip
\count255=0\global\advance\count255by\prevgraf}%
\ifnum\count255>1\box0\else
\centerline{\petit{\bf Fig.\ts#1\unskip.\
}\ignorespaces#2}\smallskip\fi}

\def\begtab#1cm#2\endtab{\par
   \ifvoid\topins\midinsert\medskip\vbox{#2\kern#1true cm}\endinsert
   \else\topinsert\vbox{#2\kern#1true cm}\endinsert\fi}
\def\begpet{\vskip6pt\bgroup\petit}
\def\endpet{\vskip6pt\egroup}
\newcount\frpages
\newcount\frpagegoal
\def\freepage#1{\global\frpagegoal=#1\relax\global\frpages=0\relax
\loop\global\advance\frpages by 1\relax
\message{Doing freepage \the\frpages\space of
\the\frpagegoal}\null\vfill\eject
\ifnum\frpagegoal>\frpages\repeat}
\newdimen\refindent
\def\begrefchapter#1{\titlea{}{\ignorespaces#1}%
\bgroup\petit
\setbox0=\hbox{1000.\enspace}\refindent=\wd0}
\def\ref{\goodbreak
\hangindent\oldparindent\hangafter=1
\noindent\ignorespaces}
\def\refno#1{\goodbreak
\hangindent\refindent\hangafter=1
\noindent\hbox to\refindent{#1\hss}\ignorespaces}
\def\endref{\goodbreak\endpet}
\def\vec#1{{\textfont1=\tams\scriptfont1=\tamss
\textfont0=\tenbf\scriptfont0=\sevenbf
\mathchoice{\hbox{$\displaystyle#1$}}{\hbox{$\textstyle#1$}}
{\hbox{$\scriptstyle#1$}}{\hbox{$\scriptscriptstyle#1$}}}}
\def\petit{\def\rm{\fam0\ninerm}%
\textfont0=\ninerm \scriptfont0=\sixrm \scriptscriptfont0=\fiverm
 \textfont1=\ninei \scriptfont1=\sixi \scriptscriptfont1=\fivei
 \textfont2=\ninesy \scriptfont2=\sixsy \scriptscriptfont2=\fivesy
 \def\it{\fam\itfam\nineit}%
 \textfont\itfam=\nineit
 \def\sl{\fam\slfam\ninesl}%
 \textfont\slfam=\ninesl
 \def\bf{\fam\bffam\ninebf}%
 \textfont\bffam=\ninebf \scriptfont\bffam=\sixbf
 \scriptscriptfont\bffam=\fivebf
 \def\sans{\fam\sansfam\ninesans}%
 \textfont\sansfam=\ninesans \scriptfont\sansfam=\sixsans
 \scriptscriptfont\sansfam=\fivesans
 \def\tt{\fam\ttfam\ninett}%
 \textfont\ttfam=\ninett
 \normalbaselineskip=11pt
 \setbox\strutbox=\hbox{\vrule height7pt depth2pt width0pt}%
 \normalbaselines\rm
\def\vec##1{{\textfont1=\tbms\scriptfont1=\tbmss
\textfont0=\ninebf\scriptfont0=\sixbf
\mathchoice{\hbox{$\displaystyle##1$}}{\hbox{$\textstyle##1$}}
{\hbox{$\scriptstyle##1$}}{\hbox{$\scriptscriptstyle##1$}}}}}

%
\let\header=Y
\let\FIG=N
\newbox\figbox
\output={\if N\header\headline={\hfil}\fi\plainoutput\global\let\header=Y
\if Y\FIG\topinsert\unvbox\figbox\endinsert\global\let\FIG=N\fi}
\let\lasttitle=N
\def\bookauthor#1{\vfill\eject
     \bgroup
     \baselineskip=22pt
     \lineskip=0pt
     \pretolerance=10000
     \authfont
     \rightskip 0pt plus 6em
     \centerpar{#1}\vskip2true cm\egroup}
\def\bookhead#1#2{\bgroup
     \baselineskip=36pt
     \lineskip=0pt
     \pretolerance=10000
     \headfont
     \rightskip 0pt plus 6em
     \centerpar{#1}\vskip1true cm
     \baselineskip=22pt
     \subhfont\centerpar{#2}\vfill
     \parindent=0pt
     \baselineskip=16pt
     \leftskip=2.2true cm
     \markfont Springer-Verlag\newline
     Berlin Heidelberg New York\newline
     London Paris Tokyo Singapore\bigskip\bigskip
     [{\it This is page III of your manuscript and will be reset by
     Springer.}]
     \egroup\let\header=N\eject}
\def\centerpar#1{{\parfillskip=0pt
\rightskip=0pt plus 1fil
\leftskip=0pt plus 1fil
\advance\leftskip by\oldparindent
\advance\rightskip by\oldparindent
\def\newline{\break}%
\noindent\ignorespaces#1\par}}
\def\part#1#2{\vfill\supereject\let\header=N
\centerline{\subhfont#1}%
\vskip75pt
     \bgroup
\textfont0=\tpfontt \scriptfont0=\tpfonts \scriptscriptfont0=\tpfontss
\textfont1=\tpmt \scriptfont1=\tbmt \scriptscriptfont1=\tams
\textfont2=\tpst \scriptfont2=\tpss \scriptscriptfont2=\tpsss
\textfont3=\tpsyt \scriptfont3=\tasys \scriptscriptfont3=\tenex
     \baselineskip=20pt
     \lineskip=0pt
     \pretolerance=10000
     \tpfontt
     \centerpar{#2}
     \vfill\eject\egroup\ignorespaces}
\newtoks\AUTHOR
\newtoks\HEAD
\catcode`\@=\active
\def\author#1{\bgroup
\baselineskip=22pt
\lineskip=0pt
\pretolerance=10000
\markfont
\centerpar{#1}\bigskip\egroup
{\def@##1{}%
\setbox0=\hbox{\petit\kern2.5true cc\ignorespaces#1\unskip}%
\ifdim\wd0>\hsize
\message{The names of the authors exceed the headline, please use a }%
\message{short form with AUTHORRUNNING}\gdef\leftheadline{%
\hbox to2.5true cc{\folio\hfil}AUTHORS suppressed due to excessive
length\hfil}%
\global\AUTHOR={AUTHORS were to long}\else
\xdef\leftheadline{\hbox to2.5true
cc{\noexpand\folio\hfil}\ignorespaces#1\hfill}%
\global\AUTHOR={\def@##1{}\ignorespaces#1\unskip}\fi
}\let\INS=E}
\def\address#1{\bgroup
\centerpar{#1}\bigskip\egroup
\catcode`\@=12
\vskip2cm\noindent\ignorespaces}
\let\INS=N%
\def@#1{\if N\INS\unskip\ $^{#1}$\else\if
E\INS\noindent$^{#1}$\let\INS=Y\ignorespaces
\else\par
\noindent$^{#1}$\ignorespaces\fi\fi}%
\catcode`\@=12
\headline={}
\def\rightheadline{}
\def\leftheadline{}

\let\header=Y
\def\contributionrunning#1{\message{Running head on right hand sides
(CONTRIBUTION)
has been changed}\gdef\rightheadline{\hfill\ignorespaces#1\unskip
\hbox to2.5true cc{\hfil\folio}}\global\HEAD={\ignorespaces#1\unskip}}
\def\authorrunning#1{\message{Running head on left hand sides (AUTHOR)
has been changed}\gdef\leftheadline{\hbox to2.5true cc{\folio
\hfil}\ignorespaces#1\hfill}\global\AUTHOR={\ignorespaces#1\unskip}}
\let\lasttitle=N
 \def\contribution#1{\vfill\supereject
 \ifodd\pageno\else\null\vfill\supereject\fi
 \let\header=N\bgroup
 \textfont0=\tafontt \scriptfont0=\tafonts \scriptscriptfont0=\tafontss
 \textfont1=\tamt \scriptfont1=\tams \scriptscriptfont1=\tams
 \textfont2=\tast \scriptfont2=\tass \scriptscriptfont2=\tasss
 \par\baselineskip=16pt
     \lineskip=16pt
     \tafontt
     \raggedright
     \pretolerance=10000
     \noindent
     \centerpar{\ignorespaces#1}%
     \vskip12pt\egroup
     \nobreak
     \parindent=0pt
     \everypar={\global\parindent=1.5em
     \global\let\lasttitle=N\global\everypar={}}%
     \global\let\lasttitle=A%
     \setbox0=\hbox{\petit\def\newline{ }\def\fonote##1{}\kern2.5true
     cc\ignorespaces#1}\ifdim\wd0>\hsize
     \message{Your CONTRIBUTIONtitle exceeds the headline,
please use a short form
with CONTRIBUTIONRUNNING}\gdef\rightheadline{\hfil CONTRIBUTION title
suppressed due to excessive length\hbox to2.5true cc{\hfil\folio}}%
\global\HEAD={HEAD was to long}\else
\gdef\rightheadline{\hfill\ignorespaces#1\unskip\hbox to2.5true
cc{\hfil\folio}}\global\HEAD={\ignorespaces#1\unskip}\fi
\catcode`\@=\active
     \ignorespaces}
 \def\contributionnext#1{\vfill\supereject
 \let\header=N\bgroup
 \textfont0=\tafontt \scriptfont0=\tafonts \scriptscriptfont0=\tafontss
 \textfont1=\tamt \scriptfont1=\tams \scriptscriptfont1=\tams
 \textfont2=\tast \scriptfont2=\tass \scriptscriptfont2=\tasss
 \par\baselineskip=16pt
     \lineskip=16pt
     \tafontt
     \raggedright
     \pretolerance=10000
     \noindent
     \centerpar{\ignorespaces#1}%
     \vskip12pt\egroup
     \nobreak
     \parindent=0pt
     \everypar={\global\parindent=1.5em
     \global\let\lasttitle=N\global\everypar={}}%
     \global\let\lasttitle=A%
     \setbox0=\hbox{\petit\def\newline{ }\def\fonote##1{}\kern2.5true
     cc\ignorespaces#1}\ifdim\wd0>\hsize
     \message{Your CONTRIBUTIONtitle exceeds the headline,
please use a short form
with CONTRIBUTIONRUNNING}\gdef\rightheadline{\hfil CONTRIBUTION title
suppressed due to excessive length\hbox to2.5true cc{\hfil\folio}}%
\global\HEAD={HEAD was to long}\else
\gdef\rightheadline{\hfill\ignorespaces#1\unskip\hbox to2.5true
cc{\hfil\folio}}\global\HEAD={\ignorespaces#1\unskip}\fi
\catcode`\@=\active
     \ignorespaces}
\def\motto#1#2{\bgroup\petit\leftskip=6.5true cm\noindent\ignorespaces#1
\if!#2!\else\medskip\noindent\it\ignorespaces#2\fi\bigskip\egroup
\let\lasttitle=M
\parindent=0pt
\everypar={\global\parindent=\oldparindent
\global\let\lasttitle=N\global\everypar={}}%
\global\let\lasttitle=M%
\ignorespaces}
\def\abstract#1{\bgroup\petit\noindent
{\bf Abstract: }\ignorespaces#1\vskip28pt\egroup
\let\lasttitle=N
\parindent=0pt
\everypar={\global\parindent=\oldparindent
\global\let\lasttitle=N\global\everypar={}}%
\ignorespaces}
\def\titlea#1#2{\if N\lasttitle\else\vskip-28pt
     \fi
     \vskip18pt plus 4pt minus4pt
     \bgroup
\textfont0=\tbfontt \scriptfont0=\tbfonts \scriptscriptfont0=\tbfontss
\textfont1=\tbmt \scriptfont1=\tbms \scriptscriptfont1=\tbmss
\textfont2=\tbst \scriptfont2=\tbss \scriptscriptfont2=\tbsss
\textfont3=\tasys \scriptfont3=\tenex \scriptscriptfont3=\tenex
     \baselineskip=16pt
     \lineskip=0pt
     \pretolerance=10000
     \noindent
     \tbfontt
     \rightskip 0pt plus 6em
     \setbox0=\vbox{\vskip23pt\def\fonote##1{}%
     \noindent
     \if!#1!\ignorespaces#2
     \else\setbox0=\hbox{\ignorespaces#1\unskip\ }\hangindent=\wd0
     \hangafter=1\box0\ignorespaces#2\fi
     \vskip18pt}%
     \dimen0=\pagetotal\advance\dimen0 by-\pageshrink
     \ifdim\dimen0<\pagegoal
     \dimen0=\ht0\advance\dimen0 by\dp0\advance\dimen0 by
     3\normalbaselineskip
     \advance\dimen0 by\pagetotal
     \ifdim\dimen0>\pagegoal\eject\fi\fi
     \noindent
     \if!#1!\ignorespaces#2
     \else\setbox0=\hbox{\ignorespaces#1\unskip\ }\hangindent=\wd0
     \hangafter=1\box0\ignorespaces#2\fi
     \vskip18pt plus4pt minus4pt\egroup
     \nobreak
     \parindent=0pt
     \everypar={\global\parindent=\oldparindent
     \global\let\lasttitle=N\global\everypar={}}%
     \global\let\lasttitle=A%
     \ignorespaces}
 \def\titleb#1#2{\if N\lasttitle\else\vskip-28pt
     \fi
     \vskip18pt plus 4pt minus4pt
     \bgroup
\textfont0=\tenbf \scriptfont0=\sevenbf \scriptscriptfont0=\fivebf
\textfont1=\tams \scriptfont1=\tamss \scriptscriptfont1=\tbmss
     \lineskip=0pt
     \pretolerance=10000
     \noindent
     \bf
     \rightskip 0pt plus 6em
     \setbox0=\vbox{\vskip23pt\def\fonote##1{}%
     \noindent
     \if!#1!\ignorespaces#2
     \else\setbox0=\hbox{\ignorespaces#1\unskip\enspace}\hangindent=\wd0
     \hangafter=1\box0\ignorespaces#2\fi
     \vskip10pt}%
     \dimen0=\pagetotal\advance\dimen0 by-\pageshrink
     \ifdim\dimen0<\pagegoal
     \dimen0=\ht0\advance\dimen0 by\dp0\advance\dimen0 by
     3\normalbaselineskip
     \advance\dimen0 by\pagetotal
     \ifdim\dimen0>\pagegoal\eject\fi\fi
     \noindent
     \if!#1!\ignorespaces#2
     \else\setbox0=\hbox{\ignorespaces#1\unskip\enspace}\hangindent=\wd0
     \hangafter=1\box0\ignorespaces#2\fi
     \vskip8pt plus4pt minus4pt\egroup
     \nobreak
     \parindent=0pt
     \everypar={\global\parindent=\oldparindent
     \global\let\lasttitle=N\global\everypar={}}%
     \global\let\lasttitle=B%
     \ignorespaces}
 \def\titlec#1#2{\if N\lasttitle\else\vskip-23pt
     \fi
     \vskip18pt plus 4pt minus4pt
     \bgroup
\textfont0=\tenbfne \scriptfont0=\sevenbf \scriptscriptfont0=\fivebf
\textfont1=\tams \scriptfont1=\tamss \scriptscriptfont1=\tbmss
     \tenbfne
     \lineskip=0pt
     \pretolerance=10000
     \noindent
     \rightskip 0pt plus 6em
     \setbox0=\vbox{\vskip23pt\def\fonote##1{}%
     \noindent
     \if!#1!\ignorespaces#2
     \else\setbox0=\hbox{\ignorespaces#1\unskip\enspace}\hangindent=\wd0
     \hangafter=1\box0\ignorespaces#2\fi
     \vskip6pt}%
     \dimen0=\pagetotal\advance\dimen0 by-\pageshrink
     \ifdim\dimen0<\pagegoal
     \dimen0=\ht0\advance\dimen0 by\dp0\advance\dimen0 by
     2\normalbaselineskip
     \advance\dimen0 by\pagetotal
     \ifdim\dimen0>\pagegoal\eject\fi\fi
     \noindent
     \if!#1!\ignorespaces#2
     \else\setbox0=\hbox{\ignorespaces#1\unskip\enspace}\hangindent=\wd0
     \hangafter=1\box0\ignorespaces#2\fi
     \vskip6pt plus4pt minus4pt\egroup
     \nobreak
     \parindent=0pt
     \everypar={\global\parindent=\oldparindent
     \global\let\lasttitle=N\global\everypar={}}%
     \global\let\lasttitle=C%
     \ignorespaces}
 \def\titled#1{\if N\lasttitle\else\vskip-\baselineskip
     \fi
     \vskip12pt plus 4pt minus 4pt
     \bgroup
\textfont1=\tams \scriptfont1=\tamss \scriptscriptfont1=\tbmss
     \bf
     \noindent
     \ignorespaces#1\ \ignorespaces\egroup
     \ignorespaces}
\let\ts=\thinspace
\def\footnoterule{\kern-3pt\hrule width 2true cm\kern2.6pt}
\newcount\footcount \footcount=0
\def\advftncnt{\advance\footcount by1\global\footcount=\footcount}
\def\fonote#1{\advftncnt$^{\the\footcount}$\begingroup\petit
\parfillskip=0pt plus 1fil
\def\textindent##1{\hangindent0.5\oldparindent\noindent\hbox
to0.5\oldparindent{##1\hss}\ignorespaces}%
\vfootnote{$^{\the\footcount}$}{#1\vskip-9.69pt}\endgroup}
\def\item#1{\par\noindent
\hangindent6.5 mm\hangafter=0
\llap{#1\enspace}\ignorespaces}

\def\titleao#1{\vfill\supereject
\ifodd\pageno\else\null\vfill\supereject\fi
\let\header=N
     \bgroup
\textfont0=\tafontt \scriptfont0=\tafonts \scriptscriptfont0=\tafontss
\textfont1=\tamt \scriptfont1=\tams \scriptscriptfont1=\tamss
\textfont2=\tast \scriptfont2=\tass \scriptscriptfont2=\tasss
\textfont3=\tasyt \scriptfont3=\tasys \scriptscriptfont3=\tenex
     \baselineskip=18pt
     \lineskip=0pt
     \pretolerance=10000
     \tafontt
     \centerpar{#1}%
     \vskip75pt\egroup
     \nobreak
     \parindent=0pt
     \everypar={\global\parindent=\oldparindent
     \global\let\lasttitle=N\global\everypar={}}%
     \global\let\lasttitle=A%
     \ignorespaces}






\def\leaderfill{\kern0.5em\leaders\hbox to 0.5em{\hss.\hss}\hfill\kern
0.5em}
\newdimen\chapindent
\newdimen\sectindent
\newdimen\subsecindent
\newdimen\thousand
\setbox0=\hbox{\bf 10. }\chapindent=\wd0
\setbox0=\hbox{10.10 }\sectindent=\wd0
\setbox0=\hbox{10.10.1 }\subsecindent=\wd0
\setbox0=\hbox{\enspace 100}\thousand=\wd0
\def\contpart#1#2{\medskip\noindent
\vbox{\kern10pt\leftline{\textfont1=\tams
\scriptfont1=\tamss\scriptscriptfont1=\tbmss\bf
\advance\chapindent by\sectindent
\hbox to\chapindent{\ignorespaces#1\hss}\ignorespaces#2}\kern8pt}%
\let\lasttitle=Y\par}
\def\contcontribution#1#2{\if N\lasttitle\bigskip\fi
\let\lasttitle=N\line{{\textfont1=\tams
\scriptfont1=\tamss\scriptscriptfont1=\tbmss\bf#1}%
\if!#2!\hfill\else\leaderfill\hbox to\thousand{\hss#2}\fi}\par}
\def\conttitlea#1#2#3{\line{\hbox to
\chapindent{\strut\bf#1\hss}{\textfont1=\tams
\scriptfont1=\tamss\scriptscriptfont1=\tbmss\bf#2}%
\if!#3!\hfill\else\leaderfill\hbox to\thousand{\hss#3}\fi}\par}
\def\conttitleb#1#2#3{\line{\kern\chapindent\hbox
to\sectindent{\strut#1\hss}{#2}%
\if!#3!\hfill\else\leaderfill\hbox to\thousand{\hss#3}\fi}\par}
\def\conttitlec#1#2#3{\line{\kern\chapindent\kern\sectindent
\hbox to\subsecindent{\strut#1\hss}{#2}%
\if!#3!\hfill\else\leaderfill\hbox to\thousand{\hss#3}\fi}\par}
\long\def\lemma#1#2{\removelastskip\vskip\baselineskip\noindent{\tenbfne
Lemma\if!#1!\else\ #1\fi\ \ }{\it\ignorespaces#2}\vskip\baselineskip}
\long\def\proposition#1#2{\removelastskip\vskip\baselineskip\noindent{\tenbfne
Proposition\if!#1!\else\ #1\fi\ \ }{\it\ignorespaces#2}\vskip\baselineskip}
\long\def\theorem#1#2{\removelastskip\vskip\baselineskip\noindent{\tenbfne
Theorem\if!#1!\else\ #1\fi\ \ }{\it\ignorespaces#2}\vskip\baselineskip}
\long\def\corollary#1#2{\removelastskip\vskip\baselineskip\noindent{\tenbfne
Corollary\if!#1!\else\ #1\fi\ \ }{\it\ignorespaces#2}\vskip\baselineskip}
\long\def\example#1#2{\removelastskip\vskip\baselineskip\noindent{\tenbfne
Example\if!#1!\else\ #1\fi\ \ }\ignorespaces#2\vskip\baselineskip}
\long\def\exercise#1#2{\removelastskip\vskip\baselineskip\noindent{\tenbfne
Exercise\if!#1!\else\ #1\fi\ \ }\ignorespaces#2\vskip\baselineskip}
\long\def\problem#1#2{\removelastskip\vskip\baselineskip\noindent{\tenbfne
Problem\if!#1!\else\ #1\fi\ \ }\ignorespaces#2\vskip\baselineskip}
\long\def\solution#1#2{\removelastskip\vskip\baselineskip\noindent{\tenbfne
Solution\if!#1!\else\ #1\fi\ \ }\ignorespaces#2\vskip\baselineskip}
\long\def\proof{\removelastskip\vskip\baselineskip\noindent{\it
Proof.\quad}\ignorespaces}
\long\def\remark#1{\removelastskip\vskip\baselineskip\noindent{\it
Remark.\quad}\ignorespaces#1\vskip\baselineskip}
\long\def\definition#1#2{\removelastskip\vskip\baselineskip\noindent{\tenbfne
Definition\if!#1!\else\
#1\fi\ \ }\ignorespaces#2\vskip\baselineskip}
\long\def\example#1#2{\removelastskip\vskip\baselineskip\noindent{\tenbfne
Example\if!#1!\else\ #1\fi\ \ }\ignorespaces#2\vskip\baselineskip}
\long\def\Examples#1{\removelastskip\vskip\baselineskip\noindent{\tenbfne
Examples\if!#1!\else\ #1\fi\ \ }}

\long\def\exercise#1#2{\removelastskip\vskip\baselineskip\noindent{\tenbfne
Exercise\if!#1!\else\ #1\fi\ \ }\ignorespaces#2\vskip\baselineskip}
\long\def\problem#1#2{\removelastskip\vskip\baselineskip\noindent{\tenbfne
Problem\if!#1!\else\ #1\fi\ \ }\ignorespaces#2\vskip\baselineskip}
\long\def\solution#1#2{\removelastskip\vskip\baselineskip\noindent{\tenbfne
Solution\if!#1!\else\ #1\fi\ \ }\ignorespaces#2\vskip\baselineskip}
\long\def\proof{\removelastskip\vskip\baselineskip\noindent{${\vec 
P\vec r\vec o\vec o\vec f.}$\quad}\ignorespaces}

\long\def\remark#1#2{\removelastskip\vskip\baselineskip\noindent{\tenbfne
Remark\if!#1!\else\ #1.\fi\ \ }\ignorespaces#2\vskip\baselineskip}
\long\def\comments#1#2{\removelastskip\vskip\baselineskip\noindent{\tenbfne
Comments\if!#1!\else\ #1\fi\ \ }\ignorespaces#2\vskip\baselineskip}

\long\def\condition#1#2{\removelastskip\vskip\baselineskip\noindent{\tenbfne\if!#1!\else\ #1\fi\ \ }\ignorespaces#2\vskip\baselineskip}

\def\frame#1{\bigskip\vbox{\hrule\hbox{\vrule\kern5pt
\vbox{\kern5pt\advance\hsize by-10.8pt
\centerline{\vbox{#1}}\kern5pt}\kern5pt\vrule}\hrule}\bigskip}
\def\typeset{\petit\noindent \par}
\outer\def\byebye{\bigskip\bigskip\typeset
\footcount=1\ifx\speciali\undefined\else
\loop\smallskip\noindent special character No\number\footcount:
\csname special\romannumeral\footcount\endcsname
\advance\footcount by 1\global\footcount=\footcount
\ifnum\footcount<11\repeat\fi
\gdef\leftheadline{\hbox to2.5true cc{\folio\hfil}\ignorespaces
\the\AUTHOR\unskip: \the\HEAD\hfill}\vfill\supereject\end}

\contribution{OPTICAL HAMILTONIANS  AND SYMPLECTIC TWIST MAPS}
\author{CHRISTOPHE GOL\'E}
\authorrunning{}
\contributionrunning{}
\address{ IMS and Department of Mathematics, SUNY at Stony Brook, Stony Brook NY 11794-3651\fonote{Partially supported by an NSF Postdoctoral 
Fellowship DMS 91-07950.}}
\abstract{ This paper concentrates on optical Hamiltonian systems of $\TTn$,
i.e. those for which $\Hpp$ is a positive definite matrix, and their
relationship with symplectic twist maps. We present theorems of decomposition
 by symplectic twist maps and existence of periodic orbits for these systems.
The novelty of these results resides in the fact that no explicit asymptotic
condition is imposed on the system. We also present a theorem of suspension
by Hamiltonian systems for the class of symplectic twist map that emerges
in our study. Finally, we extend our results to manifolds of negative curvature.}

\advsection
\titlea{\secno}{Introduction}

In a previous paper [G91b], the author explained how symplectic twist maps could
be used to decompose Hamiltonian systems on the cotangent bundle
of a compact manifold $M^n$, thus 
deriving a discrete variational approach to the search of periodic orbits
for such  systems. This method can be seen as a generalization of the
so called ``method of broken geodesics'' in differential
geometry. A similar method was introduced by Marc Chaperon for Hamiltonian
systems [Ch84]. Although our method  is very similar to his, it is in fact even more akin to the original method (see e.g. [Mi69].)

 In [G91b], we put a  boundary condition on the Hamiltonian, however:
it had to equal the metric Hamiltonian $H_0(\q,\p)={1\over 2}\norm{\p}^2$
on a fixed level set $\{H_0=K\}$. It could be anything inside $\{H_0<K\}$, 
including time dependent. 
The result (see [G91b,92b]) was then the existence of at least $cl(M)$ (or $sb(M)$ if all nondegenerate) contractible periodic orbits inside $\{H_0<K\}$ 
for such  systems. Other results  for non contractible orbits were obtained 
if $M$ supports a metric with negative curvature, or for $\Tn$ (comparable to
Theorem {\bf 7.4} in this paper.)

Here we swap the boundary condition (and the compactness of $\{H_0\leq K\}$)
for a convexity condition which gives  a Hamiltonian  system with a priori no compact invariant set:
the systems we study here are {\it optical}, in the sense that the second
derivative in the fiber direction, $\Hpp$ is positive definite . 

The main result of this paper is:
\theorem{5.2}{
Let $H(\q,\p,t)=H_t(\z)$ be a twice differentiable function on $\TTn\times \RR$  
satisfying the following:
\item{(1)} $\sup\norm {\nabla^2 H_t}< K$
\item{(2)} The matrices $\Hpp(\z,t)$ are positive definite and
$C<\norm{\Hpp}<C^{-1}$ for some $C$.

Then the time 1 map of the associated Hamiltonian
flow has at least $n+1$ distinct periodic orbits of type $\m,d$, for each prime $\m,d \in \Zn\times\ZZ$,
and at least $2^n$ in the generic case when they are all non degenerate.}

(An $\m,d$--orbit is one for which the  $d^{th}$ iterate of each point of
the orbit is a
translation by $(\m,0)$ of this point, in the covering space $\R2n$ of $\TTn$.)

Earlier results on the existence and multiplicity of periodic
orbits can be found for Hamiltonian systems  or symplectic
twist maps of $\TTn$ in [BK87], [Che92], [CZ83],[Fe89], [G91a], [J91]. The two first are perturbative 
results, i.e. for systems close to integrable ones. The four latter works are 
global in that sense, and do not require the system to be optical, but instead
require some asymptotic condition on the first derivative of $H$ (or of the
time 1 map $F$). Only the two last works consider homotopically nontrivial orbits.

 Note also that, via the Legendre transformation, Theorem {\bf 5.2}
applies to Lagrangian systems whose Lagrangian function satisfies the same conditions as $H$ in our theorem (it is not hard to see that these conditions
translate under the Legendre transformation.) Hence Theorem {\bf 5.2} extends  some existing theorems for such systems (see, e.g., [MW89], Theorem {\bf 9.3}.)

We start (Sections 2 and 3) with some background on symplectic twist maps. 
In Section 4, we give the
proof of a theorem of existence and multiplicity of periodic orbits for
compositions of symplectic twist maps
with a
convexity condition (Theorem {\bf 4.3}.) A proof of such a theorem was given in [KM89].
Unfortunately, the multiplicity part of their proof is wrong. We reproduce here
their proof of existence of a minimum, and present a new proof of the multiplicity.

In Section 5, we prove Theorem {\bf 5.2}. This results derives from the former theorem on symplectic 
twist maps, and  a decomposition technique. The resulting discrete variational method is interpreted as a method of broken geodesics.

 In Section 6 we show that
a symplectic twist map with the convexity condition  can be suspended by a Hamiltonian (our proof does not force the Hamiltonian to be optical, unfortunately ), extending a result of Moser [Mo86] and a remark about it
by Bialy and Polterovitch [BP92].

In section 7, we indicate how to extend Theorems {\bf 4.3} and {\bf 5.2} to the cotangent
bundle of a manifold of negative curvature. 

It is very likely that, with a 
little care, these techniques could extend to optical Hamiltonians on the
cotangent bundle of any compact manifolds. 

The author would like to thank J.D. Meiss and F. Tangerman for useful 
conversations, and D. McDuff and J.D. Meiss for the corrections they suggested.

\advsection
\titlea{\secno }{ Symplectic Twist Maps of $\Tn \times \Rn$}
 Let $\Tn =\Rn/\Zn$ be the $n$--dimensional torus. Its cotangent bundle
$\TTn\buildrel \pi\over\rightarrow\Tn$
is trivial: $\TTn=\Tn\times \Rn$, the cartesian product of $n$ cylinders. We give it the coordinates $(\q,\p)$ in which the symplectic structure is 
$$
\Omega=d\q\wedge d\p=\sum_{k=1}^ndq_k\wedge dp_k.
$$
 As in any cotangent bundle, $\Omega$ is exact:
$\Omega=-d\lambda$, where $\lambda= \p d\q$. 

It  is useful to work in the 
covering space $\R2n=\tilde\Tn \times \Rn$ of $\TTn$,
with projection $\pr:\R2n \to \Tn\times\Rn$.

 Of course, $\pr$ is an
exact symplectic map (see Definition 2.1) , as we have $pr^*\p d\q -\p d\q=0$.

The group $\Zn$ of
 deck or covering transformations is the set of integer vector translation in $\R2n$ of the form:
$$
\Tm.(\q,\p)=(\q +\m,\p), \ \ \ \m \in \Z^n.
$$
	A  lift of a map $F:\TTn \to \TTn$ is a map $\F:\R2n\to\R2n$ such
that $\pr\circ\F=F\circ \pr$. Since $\pr$ is a local, symplectic diffeomorphism, $\F$ is symplectic if and only if $F$ is. On the other hand,
$\F$ will always be exact symplectic when it is symplectic, which is not
the case for $F$, as the example $(\q,\p)\to (\q, \p +\p_0)$ shows.

We will  fix the lift of a map $F$ once and for all, remembering 
that two lifts only differ by a composition by  some $\Tm$.

\definition{\thmno 1}{ A map $F$ of $\TTn$ is called a {\bf
symplectic twist map} if
\item{(1)} $F$ is homotopic to $Id$.
\item{(2)} $F$  is {\bf exact symplectic}: $F^*\p d\q-\p d\q= dh$ for some $h:\TTn\to\RR$.
\item{(3)} {\bf (Twist Condition)}  If $\F(\q,\p)=(\Q,\P)$ is a lift of $F$ then
the map $\p \to \Q(\q_0,\p)$ is a diffeomorphism of $\Rn$ for
all $\q_0$, and thus the map
$$
\psi : \  (\q,\p)\to(\q,\Q)
$$
is a diffeomorphism (change of coordinates) of $\R2n$.}

\comments{\thmno 1}{ 
\item{1}.The twist condition (3) implies the more familiar looking:
$$ 
\det\ \del \Q/\del \p \ne 0.
$$
It also implies that :
$$
\F^*\p d\q-\p d\q = dS(\q,\Q)  
\numberform
$$
where $ S$ is the lift of $h$ written in the $(\q,\Q)$ coordinates
: $  S=  \t h\circ\psi^{-1}$, with $\t h=h\circ pr$. Equivalently, we
can write:
$$
\eqalign{\p=&-\del_1 S(\q,\Q)\cr
 \P=&\del_2 S(\q,\Q).\cr}
$$
  
$S(\q,\Q)$ is called a {\bf generating function} for $\F$.

\item{2}.  Condition (1) is equivalent to the fact that on {\it any}
lift $\F$ of $F$:
$$
\F\circ \Tm=\Tm\circ \F,\  \hbox{\rm i.e.} \quad \F(\q+\m,\p)=\F(\q,\p)+(\m,0).
\numberform
$$}

\example{\thmno 1}{
The  family of maps 
$$
F_s(\q,\p)= \left(\q + A\left(\p -\nabla V_s(\q)\right), \p -\nabla V_s(\q)\right)
$$
where $A$ is a nondegenerate symmetric matrix, $V_s$ is a $C^2$ function on $\Tn$ is called the {\bf  standard family}.
 Usually, $V_0\equiv 0$. The generating function for $F_s$ is
given by:
$$
S_s(\q,\Q)= S_0(\q,\Q) +V_s(\q).
$$
Where
$$
S_0(\q,\Q)={1\over 2}\langle A^{-1}(\Q-\q),(\Q-\q)\rangle
$$
is the generating function of the {\bf completely integrable map}:
$$
F_0:(\q,\p)\to (\q +A\p,\p), \quad A^t=A,\ \det A\ne 0.
$$
( The term ``completely integrable'' comes from the fact that $F_0
$ conserves each torus $\p =\p_0$, on which it acts as
a rigid ``translation''. )

This general standard family includes the classical standard family
of monotone twist maps of the annulus where $A=1$ and 
$$
V_s(q)={s\over 4\pi^2}\cos (2\pi q)
$$
 and also the Froeschl\'e family on $\TT^2\times \RR^2$ with $A=Id$ and
$$
V_s(q_1,q_2)={1\over (2\pi)^2}\{K_1cos(2\pi q_1)+K_2cos(2\pi q_2) +\lambda
cos2\pi(q_1+q_2)\} .
$$
In this case the parameter $s=(K_1,K_2,\lambda) \in \RR^3$. (See e.g. [KM89].)
 The standard map can be interpolated by an optical Hamiltonian (see Section
7).}

As this paper will suggest, many examples of symplectic twist maps
can be derived from Hamiltonian systems, and used to understand these systems.

The following results are also helpful to construct
symplectic twist maps. Their proofs can be found in [G93] (see also [H89] for
Corollary  {\bf 2.7}.)

\proposition{\thmno 1}{There is a homeomorphism between the set of lifts $\F$
of $C^1$ symplectic
twist maps of $\TTn$ and the set of $C^2$ real valued functions $S$ on $\R2n$
satisfying the following:
\item{(a)} $S(\q+\m,\Q+\m)=S(\q,\Q),\quad \forall \m \in \Zn$
\item{(b)} The maps: $\q\to \del_2 S(\q,\Q_0)$ and $\Q \to \del_1 S(\q_0,\Q)$
are diffeomorphisms of $\Rn$ for any $\Q_0$ and $\q_0$ respectively.
\item{(c)} $S(0,0)=0$.

This correspondence is given by:
$$
\F(\q,\p)=(\Q,\P)\Leftrightarrow 
\cases{\p=&$-\del_1 S(\q,\Q)$\cr \P=&$\del_2 S(\q,\Q)$.\cr}
\numberform
$$}
\lemma{\thmno 1}{ Let $f:\RR^N\to \RR^N$ be a local diffeomorphism at each point, such
that:
$$
\sup_{x\in\RR^N}\norm{(Df_x)^{-1}}=K<\infty .
$$
Then $f$ is a diffeomorphism of $\RR^N$.}
\corollary{\thmno 1}{Let $S:\quad \R2n\to \RR$ be a $C^2$ function satisfying:
$$
\eqalign{&\det \del_{12}S \ne  0\cr
&\sup_{(\q,\Q)\in \R2n}\norm {(\del_{12} S(\q,\Q))^{-1}}=K < \infty .\cr}
\numberform
$$
Then  the maps: $\q\to \del_2 S(\q,\Q_0)$ and $\Q \to \del_1 S(\q_0,\Q)$
are diffeomorphisms of $\Rn$ for any $\Q_0$ and $\q_0$ respectively, and
thus $S$ generates an exact symplectic map of $\R2n$.}

Thus , if $S$  satisfies \formno 0, as well as the periodicity condition {\it (a)} of Proposition
 \thmno {-2}, it generates a symplectic twist map.

\corollary{\thmno 1}{Let $F$ be an exact symplectic map of $\TTn$, homotopic to
$Id$. Let $\F(\q,\p)=(\Q,\P)$ be a lift of $F$.
Suppose that 
$$\sup_{\zz\in \R2n}\norm{\left(\del \Q\over\del\p\right)_\zz^{-1}}<\infty.
\numberform
$$
Then $\F$ is a symplectic twist map.}

\advsection
\titlea{\secno}{The Variational Setting}

As in the classical case of twist map ($n=1$), the generating function of a symplectic twist map
is the key to the variational setting that these maps induce. 

\proposition{\thmno 1 (Critical Action Principle)}{ Let  $F_1,\ldots,F_N$ be symplectic twist maps of $\TTn$, and let $\F_k$ be a lift of $F_k$, with generating function
$S_k$.
The sequence $\{(\q_k,\p_k)\}_{k\in\ZZ}$ is an orbit under the successive $\F_k$'s (i.e.
$\{(\q_{k+1},\p_{k+1})=\F_k(\q_k,\p_k)\}_{k\in \Z}$
, with $\F_{k+N}=\F_k,\ S_{k+N}=S_k$) if and only if the sequence $\{\q_k\}_{k\in \ZZ}$ in $\RnZ$ satisfies:
$$
\del_1S_k(\q_k,\q_{k+1}) +\del_2S_{k-1}(\q_{k-1},\q_k)=0,\quad \forall k\in \ZZ.
\numberform
$$
The correspondence is given by: $\p_k=-\del_1S_k(\q_k,\q_{k+1})$.}

 Equation \formno 0  can be interpreted formally as:
$$
\eqalign{\nabla W(\qq)&=0  \quad\quad \hbox{\rm with}\cr
W(\qq)&=\sum_{-\infty}^{\infty}\gf k{k+1}  .\cr}
$$

This interpretation makes mathematical sense when one is concerned with periodic orbits of a symplectic twist map $F$:

\definition{\thmno 1}{ A point $(\q,\p)\in \R2n$ is called a $\m,d$--{\bf point} for the 
lift $\F$ of $F$ if $\F^d(\q,\p)=(\q+\m,\p)$, where $\m \in \Zn$ and $d\in \Z$.}

Let $\F=\F_N\circ\ldots\circ\F_1$. The appropriate space of 
sequences in which to look for critical points corresponding to $\m,d$--points
of $\F$ is:
$$
\XmdN = \{\qq \in \RnZ \mid \q_{k+dN}=\q_k + \m \} 
$$
which is isomorphic to $(\RR^n)^{dN}$: the terms $(\q_1,\ldots,\q_{dN})$ determine
a whole sequence in $\XmdN$.

To find a sequence satisfying \formno 0 in $\XmdN$ is equivalent to finding
$\qq=(\q_1,\ldots,\q_{dN})$ which is a critical point for the function:
$$
\WmdN(\qq) =\sum_{k=1}^{dN}\gf k{k+1} , 
$$
in which we set $ \q_{dN+1}=\q_1+\m$. 
To see this, write :
$$
\eqalign{\p_k=-\del_1\gf k{k+1},\cr
	 \P_k=\del_2\gf k{k+1}.\cr}
$$
Then $F_k(\q_k,\p_k)=(\Q_k,\p_k)$ and 
with this notation, the proof of Proposition \thmno {-1} (for $\m,d$--points)
reduces to the suggestive:
$$
\nabla \WmdN(\qq)=\sum_{k=1}^{dN}(\P_{k-1}-\p_k)d\q_k.
$$
A little more care must be taken in order to let the topology of $\Tn$ play a
role. Note that because of the periodicity of $S$ ({\it (a)} in Proposition
{\bf 2.4}), $\WmdN$ is invariant under the $\Zn$ action on $\XmdN$:
$$
\Tm (\q_1,\ldots,\q_{dN})=(\q_1+\m,\ldots,\q_{dN}+\m).
$$

Moreover, if we want our variational approach to count $\m,d$--orbits, and
not the individual $\m,d$--points in each orbit, we should use the fact that $\WmdN$ is also invariant
under the $N$--shift map:
$$
\sigma \{\q_k\}=\{\q_{k+N}\}.
$$

Let :
$$
\xx =\XmdN/\sigma,\tau
$$ 
be the quotient of $\XmdN$ by
these two actions. We continue to call $W$ the function induced by $W$ on the
quotient $\xx$.
 
One can show ([BK87], Proposition 1 or [G91a]) that $\xx$  is the total space of a fiber bundle over $\Tn$, and
that the projection map $\XmdN\to\xx$ is a covering map.
(One makes the change of variables:
$$
\eqalign{\vv&={1\over d}\sum_1^{dN}\q_k\cr
	\vec t_k&=\q_k-\q_{k-1}-\m/d\cr}
$$
in which $\vv$ is the base coordinate, $\vec t$ the fiber.)
 In particular,
each critical point of $\WmdN$ on $\xx$ corresponds to an infinite lattice of critical points of $\WmdN$ on
$\XmdN$. Whereas the original variational problem $\nabla \WmdN=0$ on $\XmdN$ would
pick up the (infinitely many) $\m,d$--{\it points} of the lift $\F$ of $F$, 
when we restrict it to $\xx$ it exactly gives $\m,d$--{\it orbits} of $F$.

\advsection
\titlea{\secno}{ Periodic Orbits and  the Convexity Condition  }

Let $\F(\q,\p)=(\Q,\P)$ be the lift of a symplectic twist map of $\TTn$, 
and $S(\q,\Q)$ its generating function. In this section, we impose
the:
\condition{\thmno 1  Convexity Condition}{ There is a positive $a$ such that:
$$
\langle \partial_{12} S (\q,\Q).\vv,\vv \rangle \leq -a \norm{\vv}^2.
$$
uniformly in $(\q,\Q)$.}
\remark{\thmno 1}{Note that:
$$
\dd \Q\p (\q,\p)= -\left(\del_{12}S(\q,\Q)\right)^{-1},
$$
as can easily be derived by implicit differentiation of $\p=-\del_1S(\q,\Q)$.
The convexity condition \thmno {-1}  thus translates to:
$$
\left\langle \left(\dd \Q\p \right)^{-1} \vv ,\vv\right\rangle \geq a\norm {\vv}^2,\quad \forall\vv\in \Rn.
\numberform
$$
uniformly in $(\q,\p)$. This means that $F$ has bounded twist.  MacKay, Meiss and Stark [MMS89] imposed this condition on
their definition of symplectic twist maps, a terminology that we have taken from  them.}
\theorem{\thmno 1}{ Let 
$F=F_N\circ\ldots\circ F_1$  be a finite composition of symplectic twist maps
$F_k$ of $\TTn$ each satisfying the 
convexity condition. Then, for each prime $(\m,d)\in \Zn\times\ZZ$, $F$ has at least $n+1$ distinct periodic orbits of type
$\m,d$. It has at least $2^n$ of them when they are all nondegenerate.}

By a prime pair $\m,d$ we mean that at least one of the components $m_k$
of $\m$ is prime with $d$.

\remark{\thmno 1}{ One can show ([G91a]) that an $\m,d$--point of $F$ is nondegenerate
if and only if the sequence $\qq$ of $\XmdN$ it corresponds to is a 
nondegenerate critical point for $\WmdN$. The hypothesis that $F$ has
only nondegenerate $\m,d$--points is thus equivalent to the one 
that $\WmdN$ is a Morse function. Furthermore, this is a generic condition
on the
space of symplectic twist maps [G92a]. Note also that an $\m,d$ point is
also an $k\m,kd$--point for all $k\in \NN$. The reason for restricting ourselves to prime $\m,d$ is that if we were to look for $k\m,kd$ orbits, we would
also find the prescribed number of them, but with no guarantee that they would
be any different from the $\m,d$ orbits already found.}

\proof
The first part of the proof, due to Kook and Meiss, [KM89] consists in 
proving that the function $\WmdN$  is proper, and hence has a 
minimum. 

The following lemma and corollary were proven in [MMS89], and [KM89].

\lemma{\thmno 1}{ Let $S$ be the generating function of a symplectic
twist map satisfying the convexity condition \thmno {-4}. Then there is 
an $\alpha$ and positive $\beta$ and $\gamma$ such that:
$$
S(\q,\Q)\geq \alpha - \beta\norm{\q-\Q} + \gamma\norm {\q-\Q}^2.
\numberform
$$}
\proof We can write:
$$
S(\q,\Q)=S(\q,\q) +\int_0^1 \del_2S(\q,\Q_s).(\Q-\q)ds,
$$
where $\Q_s= (1-s)\q + s\Q$. Applying the same process to $\del_2S$, we get:
$$
\eqalign{
S(\q,\Q)&=S(\q,\q)+ \int_0^1 \del_2S(\Q_s,\Q_s).(\Q-\q)ds \cr 
&-\int_0^1ds\int_0^1\langle\del_{12}S(\Q_r,\Q_s).(\Q-\q),(\Q-\q)\rangle dr\cr
&\geq \alpha -\beta\norm{\Q-\q} +\gamma\norm{\Q-\q}^2,\cr}
$$
where $\alpha=\min_{\Tn}S(\q,\q)$, $\beta=\max_{\Tn}\norm{\del_2S(\q,\q)}$ and
$\gamma={a\over 2}$.\qed

\corollary{\thmno 1} { For $F$ as in Theorem \thmno {-2}, there is a minimum
for $\WmdN$ (and hence an $\m,d$--point for $F$.)}

\proof Equation  \formno 0 as applied to each $S_k$ implies that $S_k$ has
a lower bound, thus $\WmdN$ does as well. We have to prove that this lower
bound is not attained at infinity, i.e., that $\WmdN$ is a proper map.

The set $\{ (\q,\Q) \in (\Rn\times\Rn)/\Zn \mid S(\q,\Q)\leq C\}$ is compact
since \formno 0 implies that  $S\leq C$ corresponds to bounded $\norm {\q-\Q}$.
Likewise the set
$$
{\cal S}=\{\qq \in \xx \mid \WmdN(\qq)\leq C\}
$$
is compact. Hence $\WmdN$ must have a minimum in the interior of ${\cal S}$,
for $C$ big enough. This point is a critical point. \qed
\remark{\thmno 1}{
We have thus found at least one $\m,d$--orbit corresponding to a minimum of $W$.
The reader should be aware that, unlike the 1 degree of freedom
case, this does not imply that the orbit is a minimum in the sense of Aubry
(see [H89].)}
We now turn to the proof of existence of at least $n+1$ distinct orbits of type $\m,d$, and $2^n$  when they are
 all nondegenerate. 

Remember that $\xx$ is
a bundle over $\Tn$ . Let $\Sigma\cong \Tn$ be its zero section. 
Let $K=\sup_{\Sigma}\WmdN(\qq)$ . Trivially, we have:
$$
\Sigma \subset \WmdN^K\buildrel {\rm def}\over =
\{\qq \in \xx \mid \WmdN \leq K\}
$$
( since $W$ is proper, for almost every $K$, $W^K$ is a compact  manifold with boundary, by Sard's Theorem.)  From this we get the commutative diagram:
$$
\matrix{H_*(\Sigma)&\buildrel {k_*}\over {\hbox to 40pt{\rightarrowfill}} 
&H_*(\xx)\cr
	i_*\searrow&&\nearrow j_*\cr
&H_*(\WmdN^K) &\cr}
\numberform
$$
where $i,j,k $ are all inclusion maps. But $k_*=Id$ since $\Sigma$ and $\xx$
have the same homotopy type. Hence $i_*$ must be injective.

If all the $\m,d$--points are nondegenerate,  $W$ is a Morse function (a generic situation ) and  by  [Mi69], \S 3, $W^K$ has the homotopy type
of a finite CW complex, with one cell of dimension $k$ for each critical point
of index $k$ in $W^K$. In particular, we have the following Morse inequalities:
$$
\#\{\hbox{\rm critical points of index } k\}\geq b_k
$$
where $b_k$ is the $k$th Betti number of $W^K$, $b_k>{n\choose k}$ in our case
since $H_*(\Tn)\hookrightarrow H_*(W^K)$. Hence there are at least $2^n$
critical points in this nondegenerate case.

If $W$ is not a Morse function, rewrite the diagram \formno 0, but in
Cohomology, reversing the arrows. Since $k^*=Id$, $j^*$ must be injective this
time. We know that the cup length $cl(X)=cl(\Tn)=n+1$. This exactly means that
there are $n$ cohomology classes $\alpha_1,\ldots,\alpha_n$ in $H^1(\xx)$
such that $\alpha_1\cup\ldots\cup\alpha_n\ne 0$. Since $j^*$ is injective,
$j^*\alpha_1\cup\ldots\cup j^*\alpha_n\ne 0$ and thus $cl(W^K)\geq n+1$. 
$W^K$ being compact, and invariant under the gradient flow, Lusternik-Schnirelman theory implies that $W$ has at least $n+1$ critical
points in $W^K$ (The proof of Theor\`eme 1 in CH.2 \S 19 of [DNF87], which is for compact manifolds without boundaries can easily be adapted to this case.) \qed
\advsection
\theoremno=0
\titlea{\secno}{Periodic Orbits for Optical Hamiltonian Systems}

\contributionrunning{}

\condition{Assumption \thmno 1} {$H(\q,\p,t)=H_t(\z)$ is a twice differentiable function on $\TTn\times \RR$ (or $T^*M\times \RR$, where $\t M=\Rn$) and 
satisfies the following:
\item{(1)} $\sup\norm {\nabla^2 H_t}< K$
\item{(2)} The matrices $\Hpp(\z,t)$ are positive definite and
$C<\norm {\Hpp}<C^{-1}$.}

\theorem{\thmno 1}{Let $H(\q,\p,t)$ be a Hamiltonian function on $\TTn\times \RR$
satisfying Assumption \thmno {-1} . Then the time 1 map $h^1$ of the associated Hamiltonian
flow has at least $n+1$ distinct periodic orbits of type $\m,d$, for each prime $\m,d$,
and $2^n$ in the generic case when they are all non degenerate.}

\proof
we can decompose the time 1
map:
$$
h^1=h^{N\over N}\circ(h^{N-1\over N})^{-1}\circ\ldots\circ
h^{k\over N}\circ (h^{k-1\over N})^{-1}\circ\ldots\circ h^{1\over N}\circ Id.
$$
and each of the maps $h^{k\over N}\circ (h^{k-1\over N})^{-1}$ is the time ${1\over N}$ of the  (extended) flow, starting at time ${k-1\over N}$, or
in other words, the time $1/N$ of the Hamiltonian $K_t= H_{t+{k-1\over N}}$.
 Proposition {\bf 5.4} shows that, for $N$ big enough, such maps
are symplectic twist and satisfy the convexity condition {\bf 4.1} . The
result follows
from Theorem {\bf 4.3}.\qed

\remark{\thmno 1}{
Remember that Hamiltonian maps on cotangent bundles are exact symplectic.
 More precisely, the time $t$ map $h^t$ of a Hamiltonian
system on $T^*M$ satisfies:
$$
(h^t)^*\p d\q -\p d\q= dS_t \qbox{where} S_t(\q,\p)=\int_{(\q,\p)}^{h^t(\q,\p)}
\p d\q - Hds,
\numberform
$$
and the path of integration is the trajectory 
$(h^s(\q,\p),s)$ of the (extended) flow. Obviously $h^t$ is 
isotopic to $Id$. The twist condition is what remains to be checked -- it 
is clearly not always satisfied.  The following proposition shows that it is,
for small $t$, under Assumption \thmno {-2}.

Before that, let us remark that the method of proof that
we are using in this section is analogous to the so called {\it method of broken 
geodesics}  [Mi69]: by \formno 0, the function $W$ that we appeal to above in our
use of Theorem {\bf 4.3} can be interpreted as:
$$
W(\qq)=\sum_k\int_{\gamma_k}\p d\q -H ds
$$
where $\gamma_k$ is the orbit of $h^t$ starting from $(\q_k,\p_k)$ at time 
${k\over N}$, and ending at $(\q_{k+1},\P_k)$ at time ${k+1\over N}$.
The broken curve whose pieces are the $\gamma_k$ projects, via the diffeomorphisms $\psi_k$ (see definition {\bf 2.1}) to a continuous, but
only piecewise differentiable curve of $\Tn$. In the case where $H$ is the
Hamiltonian corresponding to a metric, this curve is a piecewise, or broken
geodesic and $\psi_k$ is 
the exponential map.  Proposition {\bf 3.1} can then be interpreted as saying that,
among broken geodesics, the smooth ones are exactly the ones that are
critical for $W$ (See [G93] for more details.)}

The following applies without change to Hamiltonians in cotangent bundles of
Riemannian manifolds of negative curvature. It is, however, the point at which
our method breaks for the cotangent of arbitrary manifolds: symplectic
twist maps cannot be defined on {\it all} of $T^*S^2$, for instance.

\proposition{\thmno 1}{Let $h^\epsilon$ be the time $\epsilon$ of a Hamiltonian flow
for a  Hamiltonian function satisfying  Assumption \thmno {-3}.
Then, for all sufficiently small $\epsilon$, $h^\epsilon$ is a symplectic twist map of $\TTn$. Moreover, $h^\epsilon$ satisfies
the convexity condition {\bf 4.1} .}

\proof

We can work in the covering space $\R2n$ of $T^*\Tn$, to which the
flow lifts. 
The differential of $h^t$ at a point $\vec z=(\q,\p)$ is solution
of the linear (variation) equation:
$$
\dot U(t) = J \nabla^2H(h^t(\z)) U(t),\quad U(0)=Id, \quad J=\pmatrix{0&-Id\cr
Id&0\cr}
\numberform
$$
We first need a lemma that tells us that $U(\epsilon)$ is not too far from $Id$:

\lemma{\thmno 1}{Consider the linear equation:

$$
\dot U(t)=A(t)U(t),\quad U(t_0)=U_0
$$
where $\norm {A(t)}<K, \forall t$. Then :

$$
\norm {U(t)-U_0}<K \norm{U_0}|t-t_0|e^{K|t-t_0|}.
$$}
\proof Let $V(t)=U(t)-U(t_0)$, so that $V(t_0)=0$. We have:

$$
\eqalign{\dot V(t)&= A(t)\left( U(t)-U_0\right) +A(t)U_0\cr
		  &=A(t)V(t) +A(t)U_0\cr}
$$
and hence:

$$
\norm{V(t)}=\norm{V(t)-V(0)}\leq\int_{t_0}^tK\norm{V(s)}ds +|t-t_0|K\norm {U_0}
$$
For all $|t-t_0|\leq \epsilon$, we can apply Gronwall's inequality to get:

$$
\norm{V(t)}\leq \epsilon K\norm{U_0}e^{K|t-t_0|}
$$

and we get the result by setting $\epsilon=\abs{t-t_0}$.\qed

We now finish the proof of Proposition \thmno {-1} .
By Lemma \thmno {0} we can write:

$$
U(\epsilon)-Id =\int_0^\epsilon J\nabla^2 H (h^s(\z)).(Id + O_1(s)) ds
$$

where $\norm {O_1(s)}<2Ks$, for $s\leq\epsilon$ small enough. 

Let $(\q(t),\p(t))=h^t(\q,\p)=h^t(\zz)$. The matrix
 $\b_\epsilon(\z)=\del \q(\epsilon)/\del \p$, is the upper
right $n\times n$ matrix of $U(\epsilon)$. It is given by:
$$
\b_\epsilon(\z)=\int_0^\epsilon \Hpp (h^s(\z))ds +\int_0^\epsilon O_2(s)ds
\numberform
$$
where $\left|\int_0^\epsilon O_2(s)ds\right|< K\epsilon^2$. From this,
and the fact that 
$$C\norm {\vec v}^2<\langle \Hpp(\z) \vec v,\vec v\rangle <C^{-1}\norm {\vec v}^2,$$
we deduce that:
$$
 (\epsilon C -K\epsilon^2)
\norm {\vv}^2<\langle \b_\epsilon(\zz ) \vv ,\vv\rangle < (\epsilon C^{-1} +K\epsilon^2)
\norm {\vv}^2
\numberform
$$
so that in particular $\b_\epsilon(\z)$ is nondegenerate for small enough $\epsilon$.
The set of nonsingular matrices $\{\b_\epsilon(\z)\}_{\z\in \R2n}$ is included in
a compact set and thus:
$$
\sup_{\z\in \R2n}\norm{\b^{-1}_\epsilon(\z)}<K',
\numberform
$$
for some positive $K'$.
We can now apply Corollary {\bf 2.7} to show that $h^\epsilon$ 
is a symplectic twist map with a generating function $S$ defined
on all of $\R2n$.

Likewise, from \formno {-2}, and the 
fact that $\langle \Hpp^{-1}(\zz)\vv,\vv \rangle>C{\norm\vv}^2$, one easily derives that $h^\epsilon$ satisfies the convexity condition  {\bf 4.1}. 
\qed

\advsection
\titlea{\secno}{ Suspension of Symplectic Twist Maps by Hamiltonian Flows}
\contributionrunning{\it \the\sectionno. Suspension of Symplectic Twist Maps}

 In [Mo86], Moser showed how to suspend a monotone twist map of the compact
annulus into a time 1 map of a (time dependent) optical
 Hamiltonian system. Furthermore, he was careful to construct the Hamiltonian in such a way that its flow leaves invariant the compact annulus (when the map does) and also such that it is time periodic.

As announced by Bialy and Polterovitch [BP92], Moser's method can be adapted to suspend a
symplectic twist map whose generating functions $S$ is such that
$\del_{12}S(\q,\Q)$ is a  positive definite  
{\it symmetric} matrix satisfying the convexity condition {\bf 4.1}.
 In particular their result shows that
the Standard Map is the time 1 of an optical Hamiltonian flow, periodic in time.

Here we present a suspension theorem for higher dimensional symplectic
twist maps, without
the assumption that $\del_{12}S$ is
 symmetric.  Our result is modest in that we do not
obtain a convexity condition on the Hamiltonian, or show that the Hamiltonian
we construct can be made time periodic. Our method is different from Moser's.

\theorem{\thmno 1} {Let $F(\q,\p)=(\Q,\P)$ be a symplectic twist map of $\TTn$ 
which  satisfies the convexity condition {\bf 4.1}. 
Then $F$ is the time 1 map of a (time dependent) Hamiltonian $H$.}

\proof
Let $S(\q,\Q)$ be the generating function of $F$. Condition {\bf 4.1} can be
rewritten:
$$
\inf_{(\q,\Q)\in \R2n}\langle -\del_{12}S(\q,\Q) \vv,\vv\rangle > a\norm \vv^2, \quad a>0, \forall \vv\ne 0 \in \Rn.
\numberform
$$
The following lemma, whose proof is left to the reader shows that this inequality implies (2.4). Hence
whenever we have a function on $\R2n$ which is suitably periodic and satisfies
\formno 0, it is the generating function for some symplectic twist map.

\lemma{ \thmno 1 }{Let $\{A_x\}_{x\in \Lambda}$ be a family of 
 $n\times n$ real matrices satisfying:
$$
\sup_{x\in \Lambda}{\langle A_x \vv,\vv\rangle }> a{\norm \vv}^2, \quad \forall
\vv \ne 0 \in \Rn.
$$
Then :
$$
\sup_{x\in \Lambda}\norm {A_x^{-1}} < a^{-1}.
$$}

We  construct a differentiable family  $S_t$ of generating functions, with $S_1=S$,
and then show how to make a Hamiltonian vector field out of it, whose
time 1 map is $F$.
Let
$$
S_t(\q,\Q)=\cases{{1\over 2}af(t){\norm {\Q-\q}}^2 &for $0<t\leq {1\over 2}$\cr
	   {1\over 2} af(t){\norm {\Q-\q}}^2 + (1-f(t))S(\q,\Q) &for ${1\over 2}\leq t\leq 
1.$\cr}
$$
where $f$ is a smooth positive functions,  $f(1)=f'(1/2)=0, \  f(1/2)=1,\ \lim_{t\to 0^+}f(t)=+\infty$.
We will ask also that $1/f(t)$, which can be continued to $1/f(0)=0$ be
differentiable at $0$. The choice of $f$ has been made so that $S_t$ is differentiable
with respect to $t$, for $t\in (0,1]$.  Furthermore, it is easy to verify that:
$$
\sup_{(\q,\Q)\in \R2n}\langle -\del_{12}S_t(\q,\Q) \vv,\vv\rangle > a\norm \vv^2, \quad a>0, \forall \vv\ne 0 \in \Rn, t\in (0,1].
$$
Hence $S_t$ generates a smooth family $F_t, \ t\in (0,1]$ of symplectic twist maps, and
in fact $F_t(\q,\p)=(\q + (af(t))^{-1}\p,\p), \quad t\leq 1/2)$, so that
$\lim_{t\to 0^+}F_t=Id$, in any topology that one desires (on compact sets.)
Let us  write
$$
\t S_t(\q,\p)= S_t\circ \psi_t (\q,\p),
$$
 where $\psi_t$ is the change of coordinates  given by the fact that $F_t$ is twist. It is
not
hard to verify that $\psi_t (\q,\p)=(\q,\q -(af(t))^{-1}\p), \quad t\leq 1/2$.
so that:
$$
\t S_t(\q,\p)= {1\over 2}(af(t))^{-2}{ \norm \p}^2
$$
In particular, by our assumption on $1/f(t)$, $\t S_t$ can be differentiably 
continued for
all $t\in [0,1]$, with $S_0\equiv 0$.
 Hence, in the $\q,\p$ coordinates, we can write:
$$
F_t^*\p d\q -\p d\q= d\t S_t, \quad t\in [0,1].
$$
A familly of maps that satisfies this with $\t S_t$ differentiable in $(\q,\p,t)$
is called an {\it exact symplectic isotopy}.
The proof of the theorem derives from the standard:
\lemma{\thmno 1}{ Let $g_t$ be an exact symplectic isotopy of $\TTn$ (or $T^*M$,
in general.) Then $g_t$ is a Hamiltonian isotopy.}
\proof 
Let $g_t$ be an exact symplectic isotopy:
$$
g_t^*\p d\q-\p d\q=dS_t
$$
for some $S_t$ differentiable in all of $(\q,\p,t)$. We  claim that the (time dependent ) vector field:
$$
X_t(\z)={dg_t\over dt}(g_t^{-1}(\z))
$$
whose time $t$ is $g_t$, is Hamiltonian. To see this, we compute:
$$
{d\over dt}(d\t S_t)={d\over dt}g_t^*\p d\q=g_t^*L_{X_t}\p d\q=
g_t^*\left(i_{X_t}d(\p d\q)-d(i_{X_t}\p d\q)\right),
$$
from which we get
$$
i_{X_t}d\q \wedge d\p= d H_t
$$
with 
$$
H_t=\left( (g_t^{-1})^*{dS_t\over dt}-i_{X_t}\p d\q \right),
$$
 which exactly means that $X_t$ is Hamiltonian.

\qed

\advsection
\titlea{\thmno 1}{ Cotangent Bundle of  Manifolds with Negative Curvature}

We indicate in this section how some of the previous results can be obtained in  the cotangent bundle $T^*M$ of
a compact manifold $M$ which
supports a metric of negative curvature. Such a manifold  is always covered by $\Rn$ ( As before we denote by $\pr:\t M ( =\Rn)\to M$ the
covering map.) The definition of symplectic twist map carries through verbatim for the cotangent bundle of such manifolds, as well as Propositions 2.4 and
3.1, Corollaries 2.6 and 2.7. The action  by translations of $
\pi(\Tn)=\Zn$ on $\R2n$  is replaced by the more general action of
$\pi_1(M)$, the deck transformation group of $T^*\t M$. Note also that the convexity
condition {\bf 4.1} still makes sense in this more general context. For more details, 
see [G91b], [G93].

The first resistance we encounter to an extension of our results to such manifolds is Definition {\bf 3.2} of $\m,d$--orbits. The clue to define
such an orbit in this new context is Remark {\bf 5.4}: we saw there that, 
in the case where the map $F$ considered is Hamiltonian and decomposed into
symplectic twist maps, an $\m,d$--sequence gives rise to a closed, piecewise
smooth curve in $\Tn$ (a ``broken geodesic''). The integer vector $\m$ 
classifies these broken geodesics up to homotopy {\it with or without fixed base 
points}. This is because the group $\pi_1(\Tn)=\Zn$ is abelian. 

In general
manifolds, two loops through a base point that represent different elements
in $\pi_1(M)$ might be homotopic if we allow the homotopy to move the base point: we say then that the curves are {\it free homotopic}. Free homotopy
classes are in one to one correspondence with the conjugacy classes in
$\pi_1(M)$.

Coming back to our broken geodesics, the natural classification for periodic
orbits of a Hamiltonian system is that of free homotopy class: each of these
classes represent a connected component in the loop space. This  motivates:

\definition{\thmno 1}{Let $\m $ be a representative of a free homotopy class of
loops in $M$.  A sequence $\{\q_k\}$ of points in $ \t M=\Rn$ is called a $\m,d$--{\bf sequence} if, for all $k\in \ZZ$,  $\pr (\q_k)=\pr(\q_{k+d})$ 
and (any) curve $\t\gamma$ of $\t M$ that joins
$\q_k$ and $\q_{k+d}$, projects to a closed curve of $M$ in the free homotopy class $\m$, independent of $k$. 
 The orbit $\{(\q_k,\p_k)\}$ of a map of $T^*M$ is an $\m,d$--{\bf periodic 
orbit} if the sequence $\{\q_k\}$ is an $\m,d$--sequence.}

We can now state:

\theorem{\thmno 1}{ Let $M$ be a compact Riemannian manifold with negative curvature.  Let 
$F=F_N\circ\ldots\circ F_1$  be a finite composition of symplectic twist maps
$F_k$ of $T^*M$ satisfying the 
convexity condition {\bf 4.1}. Then, for each free homotopy class $\m$ and
period $d$, $F$ has at least 2 periodic orbits of type
$\m,d$. If $\m=0$, the class of contractible loops, then there are at least
$cl(M)$ orbits of type $\m,d$, and $sb(M)$ if they are all nondegenerate.}

\proof It is shown in [G91b], Lemma {\bf 7.2.2}, that the set $X$ of $\m,dN$ sequences
(modulo $\pi_1(M)$ and shift; X is denoted $O_{m,d}/\sigma$ in that paper), has a deformation retraction onto the set, that we call $\Sigma$, formed by  the unique geodesic of class $\m$ (remember that $M$ has negative curvature) that is, $\xx$ has the homotopy type of $\T1$. The proof of Theorem {\bf 4.3} can 
now be repeated, keeping
in mind that $sb(\T1)=cl(\T1)=2$. 

When $\m$ is the trivial class, the set $X$ retracts 
on the set of constant loops, naturally embedded in it
([G91b], Lemma {\bf 6.2}). This set, that we call $\Sigma$ again is
homeomorphic to $M$. A simple adaptation of Lemma {\bf 7.2.2} in [G91b] shows
 that $\Sigma$ is in fact a deformation retract of $X$ and hence once again,
 we can repeat the proof of Theorem {\bf 4.3}.\qed

Assumption {\bf 5.1} and Proposition {\bf 5.4} apply without a change to our new context
and hence we have:

\theorem{\thmno 1}{ Let $M$ be a compact Riemannian manifold with negative curvature. Let $H(\q,\p,t)$ be a Hamiltonian function on $T^*M$
satisfying Assumption {\bf 5.1} . Then the time 1 map  of the associated Hamiltonian
flow can
be decomposed into a product of symplectic twist maps. It has at least $2$ periodic orbits of type $\m,d$, for each  $\m,d$. When $\m$ is the trivial class,  there are at least $cl(M)$ orbits of type $\m,d$, and $sb(M)$ if
they are all nondegenerate (i.e. generically.)}

\begrefchapter{References}
\ref{[A78]} V.I. Arnold: ``Mathematical Methods of Classical Mechanics''
 (Appendix 
 9), Springer-Verlag 1978.

\ref{[AL83]} S. Aubry and P.Y. LeDaeron: ``The discrete Frenkel-Kontarova
model and its extensions I. Exact results for ground states'', {\it Physica} 8D
 (1983), 381-422.

\ref{[BK87]} D. Bernstein and A.B. Katok: ``Birkhoff periodic orbits for small
perturbations of completely integrable Hamiltonian systems with convex
Hamiltonians'', {\it Invent. Math.} {\bf 88} (1987), 225-241.

\ref{[BP92]} M. Bialy and L. Polterovitch,
`` Hamiltonian systems, Lagrangian tori and Birkhoff's Theorem'',
 {\it Math. Ann.}
292, (1992)   619--627.

\ref{[Cha84]} M. Chaperon, ``Une id\'ee du type ``g\'eod\'esiques bris\'ees''
pour les syst\`emes hamiltoniens'', {\it C.R. Acad. Sc., Paris}, 298, S\'erie I,
no 13, (1984) 293-296.

\ref{[Che92]} Chen, Weifeng, ``Birkhoff periodic orbits for small
perturbations of completely integrable Hamiltonian systems with nondegenerate
Hessian'', in  Twist mappings and their applications,
 R R. McGehee and K. R. Meyer ( Ed.), {\it IMA Vol. in Math. and App.}
No. 44 (1992).

\ref{[CZ83]} C.C. Conley and E. Zehnder: ``The Birkhoff-Lewis 
fixed point theorem
 and a conjecture of V.I. Arnold'', {\it Invent. Math.} (1983).

\ref{[DNF87]} B. Doubrovine, S. Novikov, A. Fomenko : ``G\'eom\'etrie
contemporaine'', vol 3, Editions Mir, Moscow, (1987) (see also the English
translation in Springer--Verlag)

\ref{[Fe89]}  P. Felmer, ``Periodic solutions of spatially periodic hamiltonian systems'',  {\it CMS Technical Summary report} No. 90--3, U. of Wisconsin 
(1989).

\ref{[G91a]} C. Gol\'e, ``Monotone maps and their periodic orbits'', in
The geometry of Hamiltonian systems, T. Ratiu, (Ed.) Springer (1991).

\ref{[G91b]} C. Gol\'e, ``Periodic orbits for Hamiltonian systems in cotangent 
bundles'', {\it IMS preprint}, SUNY at Stony Brook (1991).

\ref{[G92a]} C. Gol\'e, ``Ghost circles for twist maps'',
{\it Jour. of Diff. Eq.} , Vol. 97, No. 1 (1992),   140--173.

\ref{[G92b]} C. Gol\'e, ``Symplectic twist maps and the theorem of 
Conley--Zehnder for general cotangent bundles,'' in 
Mathematical Physics X, (proceedings of the ${\rm X}^{th}$ congress)  K. Schm\"udgen (Ed.), Springer--Verlag (1991).

\ref{[G93]} C.Gol\'e, ``Symplectic twist maps'' World Scientific Publishers
(to appear).

\ref{[He89]} M.R. Herman, ``Inegalit\'es a priori pour des tores invariants
par des diff\'eomorphismes symplectiques'', {\it Publ. Math. I.H.E.S.}, No. 70 (1989),
47--101.

\ref{[J91]} F.W. Josellis: ``Global periodic orbits for Hamiltonian systems
on $\Tn\times\Rn$''
Ph.D. Thesis Nr. 9518, ETH Z\"urich, (1991).

\ref{[KM89]} H. Kook and J. Meiss: ``Periodic orbits for Reversible,
 Symplectic
 Mappings'', {\it Physica D } 35 (1989) 65--86.

\ref{[L89]} P. LeCalvez: ``Existence d'orbites de Birkhoff G\'en\'eralis\'ees
pour les diff\'eo- morphismes conservatifs de l'anneau'', Preprint, Universit\'e
Paris-Sud, Orsay (1989).

\ref{[MMS89]} R.S. MacKay, J.D. Meiss and J.Stark, ``Converse KAM theory for
symplectic twist maps'', {\it Nonlinearity } No.2 (1989) 555--570.

\ref{[Mi69]} J. Milnor: ``Morse Theory'' (second edition), Princeton University Press
 (1969).

\ref{[Mo77]} J.Moser: ``Proof of a generalized form of a fixed point
theorem due to G.D. Birkhoff'', {\it Lecture Notes in Mathematics,} Vol. 597:
Geometry and Topology, pp. 464-494. Springer (1977).

\ref{[Mo86]} J. Moser: ``Monotone twist mappings and the calculus of
 variations'',
{\it Ergod. Th. and Dyn. Sys.}, Vol 6 (1986), 401-413.

\ref{[MW89]} J. Mawhin and M. Willem: ``Critical point theory and Hamiltonian
systems'', Springer--Verlag (1989).

\endref
\centerline{Email: gole@math.sunysb.edu}

\byebye